\documentclass[preprint]{elsarticle}
\journal{Journal of Approximation Theory}

\usepackage{amsmath,amssymb,amsthm,url,times,color}
\usepackage{mathtools}

\def\today{\number\day\ \ifcase\month\or
	January\or February\or March\or April\or May\or June\or
	July\or August\or September\or October\or November\or December\fi
	\space \number\year}
\def\gobble#1#2{}
\def\shortdate{\number\day/\number\month/\expandafter\gobble\number\year}

\def\nn{\mathbb{N}}
\def\R{\mathbb{R}}

\newcommand{\hypergeom}[5]{\mbox{$
_#1 F_#2\left( \! \left.
\begin{array}{c}
\multicolumn{1}{c}{\begin{array}{c} #3
\end{array}}\\[1mm]
\multicolumn{1}{c}{\begin{array}{c} #4
            \end{array}}\end{array}
\! \right| \displaystyle{#5}\right) $} }

\def\a{\alpha}
\def\b{\beta}

\def\la{\lambda}

\newtheorem{theorem}{Theorem}[section]

\theoremstyle{definition}

\newtheorem{remark}[theorem]{Remark}

\numberwithin{figure}{section}
\numberwithin{equation}{section}
\numberwithin{table}{section}

\def\l{\lambda}

\def\beq{\begin{equation}}
	\def\eeq{\end{equation}}

\definecolor{dkg}{rgb}{0,0.7,0}
\definecolor{dkr}{rgb}{0.9,0,0}
\definecolor{dkb}{rgb}{0,0,0.7}
\definecolor{purple}{rgb}{0.5,0,0.7}

\definecolor{gold}{rgb}{0.83, 0.69, 0.22}

\def\url#1{\href{#1}{#1}}
\begin{document}
\begin{frontmatter}
\title{Bounds for extreme zeros of Meixner-Pollaczek  polynomials}	
\author[aj]{A.S. Jooste\corref{cor1}}
\ead{alta.jooste@up.ac.za}
\address[aj]{Department of Mathematics and Applied Mathematics,\\ University of Pretoria, Pretoria, South Africa}
\author[kj]{K. Jordaan}
\ead{jordakh@unisa.ac.za}
\cortext[cor1]{Corresponding author}
\address[kj]{Department of Decision Sciences,\\
University of South Africa, Pretoria, 0003, South Africa}

\begin{abstract} 
\noindent In this paper we consider connection formulae for orthogonal polynomials in the context of Christoffel transformations for the case where a weight function, not necessarily even, is multiplied by an even function $c_{2k}(x),k\in \nn_0$, to determine new lower bounds for the largest zero and upper bounds for the smallest zero of a Meixner-Pollaczek polynomial. When $p_n$ is orthogonal with respect to a weight $w(x)$ and $g_{n-m}$ is orthogonal with respect to the weight $c_{2k}(x)w(x)$, we show that $k\in\{0,1,\dots,m\}$ is a necessary and sufficient condition for existence of a connection formula involving a polynomial $G_{m-1}$ of degree $(m-1)$, such that the $(n-1)$ zeros of $G_{m-1}g_{n-m}$ and the $n$ zeros of $p_n$ interlace. We analyse the new inner bounds for the extreme zeros of Meixner-Pollaczek polynomials to determine which bounds are the sharpest. We also briefly discuss bounds for the zeros of Pseudo-Jacobi polynomials.
\end{abstract} 


\end{frontmatter}
\section{Introduction}

A sequence of real polynomials $\{p_n\}_{n=0}^\infty$, where $p_n$ is of exact degree $n,$ is orthogonal with respect to a positive measure $\mu(x)>0$ on an interval $[a,b],$ if the scalar product satisfies
\begin{equation}\nonumber
\langle p_m,p_n\rangle=
\int_{a}^{b}p_m(x)p_n(x)d\mu(x)=0,\quad m\neq n.
\end{equation}
If $\mu(x)$ is absolutely continuous, then it can be represented by a real weight function $w(x)>0$ so that $d\mu(x)=w(x)\,dx$. Zeros of orthogonal polynomials on the real line are all real and distinct and have the interesting property that the zeros of consecutive polynomials in the sequence interlace. Stieltjes showed that between any two consecutive zeros of $p_m$, $m<n-1$, there is at least one zero of $p_n$. De Boor and Saff \cite{DBS} extended this by proving that there is a polynomial $S$ of degree $(n-m-1)$, associated naturally with $p_n$, such that the $(n-1)$ zeros of $S(x)p_m(x)$ and the $n$ zeros of $p_n(x)$ separate each other. 

The question whether Stieltjes interlacing occurs between the zeros of two polynomials $p_n$ and $q_m$, $m<n-1$, from different orthogonal sequences $\{p_n\}_{n=0}^{\infty}$ and $\{q_n\}_{n=0}^\infty$ and whether polynomials analogous to the de Boor-Saff polynomials exist in this more general situation, was answered in the affirmative for different sequences $C_{n+1}^{(\lambda)}$ and $C_{n-2}^{(\lambda+k)}$, $k\in\{0,1,2,3\}$, of Gegenbauer polynomials in \cite{DriverNumerMath} and $L_{n+1}^{(\alpha)}$ and $L_{n-1}^{(\alpha+t)}$, $\alpha>-1$, $t\in\{0,1,2,3,4\}$, of Laguerre polynomials in \cite{DJLag}. 

In \cite[Thm 3]{DriverNumerMath}, Driver showed that 
for a Gegenbauer polynomial $C_n^{(\lambda)}$, $\lambda>-\tfrac 12$, $\lambda \neq 0$ of degree $n$ and for each $k\in\{1,2,\dots,n-1\}$, there exist polynomials $g_k$, $G_k$ and $H_{k-1}$ such that \[g_k(x)C_{n-k}^{(\lambda+k)}(x)=G_k(x)C_n^{(\lambda)}(x)+H_{k-1}(x)C_{n+1}^{(\lambda)}(x)\] where $g_k(x)=2^k\lambda(\lambda+1)\dots(\lambda+k-1)(1-x^2)^k$ and $\deg(G_k)=\deg(H_k)=k$. In addition, if $C_{n+1}^{(\lambda)}$ and $C_{n-k}^{(\lambda+k)}$ have no common zeros then, for each fixed $k\in\{1,\dots,n-1\}$, the $n$ real zeros of the product $C_{n-k}^{(\lambda+k)}G_k$ interlace with the $n+1$ zeros of $C_{n+1}^{(\lambda)}$.

The authors of \cite{DJLag} noted that the zeros of the de Boor-Saff polynomial analogues, that completed the Stieltjes interlacing of Laguerre polynomials from different sequences when $t\in\{3,4\}$, were the same as known upper bounds for the smallest zero of Laguerre polynomials obtained by Gupta and Muldoon in \cite[(2.9) and (2.10)]{Gupta_Muldoon_2007}. This observation, that the zeros of de Boor-Saff analogues appearing in the mixed three-term recurrence relations satisfied by classical orthogonal polynomials corresponding to different parameters are inner bounds for the extreme zeros was explicitly stated  as a corollary in \cite{DJ}. In the general setting, where $\{p_n\}_{n=0}^{\infty}$ is any sequence of orthogonal polynomials, if a polynomial $g_{n-m}(x)$, of degree $n-m$, where $m<n$, $m$, $n$ fixed, satisfies a mixed recurrence equation of the form \begin{equation}\label{mrr}f(x)g_{n-m}(x)=G(x)p_{n-1}+H(x)p_{n}(x),\end{equation} with $f(x)\neq0$, deg$(G(x))=m-1$, then the largest (smallest) zero of $G(x)$ is a strict lower (upper) bound for the largest (smallest) zero of $p_n(x)$ (cf. \cite[Thm 2.1 and Cor. 2.2]{DJ}).

The mixed recurrence equations of type \eqref{mrr} that are used to obtain inner bounds for extreme zeros using the approach outlined above, typically arise from the contiguous relations satisfied by hypergeometric functions \cite{rain}. For orthogonal polynomials of hypergeometric type in the Askey and $q$-Askey scheme (cf. \cite{KLS}), these contiguous relations  translate into mixed recurrence equations between parameter shifted orthogonal polynomials. The mixed recurrence equations can also be viewed as formulae connecting different sequences of orthogonal polynomials in the framework of Christoffel transformations (cf. \cite{Dimitrov}). Christoffel \cite{Christoffel} proved this important connection between polynomials orthogonal with respect to a weight and polynomials orthogonal with respect to a modification of the weight where the modification involves multiplication of the weight function by a polynomial. A generalisation of Christoffel's formula is proved by Szeg\H{o} in \cite[Thm 2.5]{szego} where he also discusses a connection formula for the special case when an even weight is modified by an even factor. Christoffel's theorem can be difficult to implement computationally and Gautschi \cite{Gautschi} describes a constructive algorithmic approach. 

Historically, results on zeros of classical continuous and discrete orthogonal polynomials focused more on outer bounds for extreme zeros, see for example  \cite{Area_et_al_2004, Area_et_al_2013, Dimitrov_Nikolov_2010, Dimitrov_Rafaeli_2009, Ismail_Li_1992, JT, Krasikov_2002, Krasikov_2006, Krasikov_Zarkh_2009,  szego}, but some results on inner bounds for extreme zeros of  Gegenbauer, Laguerre  and  Jacobi polynomials (cf. \cite{Area_et_al_2012, Gupta_Muldoon_2007, Krasikov_2006, szego}) and for the discrete orthogonal Charlier, Meixner, Krawtchouk and Hahn polynomials (cf. \cite{Area_et_al_2013, Krasikov_Zarkh_2009}) were known. The approach based on equations of type \eqref{mrr} described above, provided a new point of view and has been used to obtain new inner bounds for extreme zeros of Jacobi, Laguerre and Gegenbauer polynomials  \cite{DJ}, Meixner and Krawtchouk polynomials \cite{Jooste_Jordaan}, Pseudo-Jacobi polynomials \cite{JT}, Pseudo-Ultraspherical \cite{DM1}, Hahn, Continuous Hahn and Continuous Dual Hahn polynomials \cite{JNK_2017} and little $q$-Jacobi and $q$-Laguerre in \cite{Swami,Kar}. 

In this paper we determine new inner bounds for the extreme zeros of Meixner-Pollaczek polynomials and also briefly discuss known bounds for the extreme zeros of Pseudo-Jacobi polynomials. Meixner-Pollaczek, Pseudo-Jacobi (and also Continuous Dual Hahn and Wilson) polynomials are examples of polynomial sequences where the polynomial $g_{n-m,k}$ in \eqref{mrr}, obtained from making a parameter shift of $k$ units, $k\in\{0,1,\dots,m\}$, is orthogonal with respect to a polynomial of degree $2k$ times the weight function of $p_n$. Hence, for these classes of polynomials, the parameter shifted polynomials required to obtain inner bounds for the extreme zeros in the equations of type \eqref{mrr} involve an even polynomial perturbation of the orthogonality measure in the context of Christoffel transformations. 

We will begin our investigation in \S \ref{Stieltjes} by modifying \cite[Thm 3]{ASJ},  which dealt with the case where the parameter shift is equal to the degree of the perturbation, to determine necessary and sufficient conditions for the existence of mixed recurrence equations  for the case where the degree of the perturbation is double the parameter shift.  This will allow us to determine, in \S \ref{MP1}, all equations of type \eqref{mrr} that exist for Meixner-Pollaczek polynomials when the degree difference is fixed. In \S \ref{innerbounds}, we determine all the inner bounds for extreme zeros of these polynomials arising from the mixed recurrence equation technique using degree differences of two and three. We will (analytically and numerically) compare the new bounds obtained to each other, and also to outer bounds previously obtained (cf. \cite{Ismail_Li_1992}), in order to evaluate the sharpness of the new inner bounds. Finally, we discuss bounds for extreme zeros of Pseudo-Jacobi polynomials in \S \ref{PJ} where we recover some inner bounds obtained in \cite{JT} and  correct a typographical error for another. 
\section{Existence of connection formulae of a specific type}\label{Stieltjes}
In this section, we consider the sequence of polynomials $\{p_n\}_{n=0}^{\infty}$ orthogonal with respect to a weight $w(x)$, and the sequence of polynomials $\{g_{n,k}\}_{n=0}^{\infty}$, orthogonal with respect to a weight $c_{2k}(x)w(x)$, where $c_{2k}$ is an even polynomial of degree $2k$, $k\in\nn_0$ and $m\in \{2,3,\dots,n\}$. 

It is well-known that, if $\{p_n\}_{n=0}^{\infty}$ is a sequence of polynomials, satisfying the three-term recurrence equation
\begin{equation}\label{3trr}{p_{n}(x)=(x-C_{n})p_{n-1}(x)-\l_{n}p_{n-2}(x)},\end{equation}
with
$$C_n=\frac{\langle xp_{n-1},p_{n-1}\rangle}{\langle p_{n-1},p_{n-1}\rangle}~~\mbox{and}~~ \l_n=\frac{\langle xp_{n-1},p_{n-2}\rangle}{\langle p_{n-2},p_{n-2}\rangle},$$
 then, given $n$, there exists a sequence of real (associated) orthogonal polynomials $S_m^{(n)}(x), m\in\{0,1,2,\dots,n\},$ of exact degree $m$. These associated polynomials, also known as dual polynomials (cf. \cite{chiharabook,DBS,VZ}), are completely determined by the coefficients in the three-term recurrence equation and are also considered in the context of Christoffel transformations in \cite{Marcellan}. The associated polynomials satisfy the three-term recurrence equation (cf. \cite[(13)]{Beardon} and \cite[Lemma 1]{ASJ})
\begin{equation}\label{a3trr}S_{m}^{(n)}(x)=
(x-C_{n-(m-1)})S_{m-1}^{(n)}(x)-\l_{n-(m-2)}S_{m-2}^{(n)}(x),m\in\{1,2,\dots,n\},\end{equation}
 with $S_0^{(n)}(x)=1$ and $S_{-1}^{(n)}(x)=0$. 
 
 Starting with (\ref{3trr}), then lowering the index $n$, and iterating, we obtain
\begin{eqnarray}
\l_np_{n-2}(x)&=&\nonumber(x-C_n)p_{n-1}(x)-p_{n}(x)\\
\l_{n-1}p_{n-3}(x)&=&\nonumber(x-C_{n-1})p_{n-2}(x)-p_{n-1}(x)\\
\l_n\l_{n-1}p_{n-3}(x)&=&\nonumber(x-C_{n-1})\l_n p_{n-2}(x)-\l_n p_{n-1}(x)\\
&=&\nonumber(x-C_{n-1}){\Big((x-C_n)p_{n-1}(x)-p_{n}(x)\Big)}-\l_np_{n-1}(x)\\
&=&\nonumber{\Big((x-C_{n-1})(x-C_n)-\l_n \Big)}p_{n-1}(x)-{(x-C_{n-1})}p_{n}(x)\\
&=&\nonumber {S_2^{(n)}(x)}p_{n-1}(x)-S_1^{(n-1)}(x)p_{n}(x)
\end{eqnarray}
such that, for $m\in\{2,3,\dots,n\}$,
\begin{equation}\l_{n}\l_{n-1}\ldots \l_{n-m+2}p_{n-m}(x)=S_{m-1}^{(n)}(x)p_{n-1}(x)-S_{m-2}^{(n-1)}(x)\label{Beardon}p_{n}(x).\end{equation}
Furthermore, it follows from \eqref{a3trr} that $S_1^{(k)}(x)=x-C_k$, $k$ an integer, 
and this, together with iteration of \eqref{3trr} with $n$ replaced by $n+1$, yields the recurrence equation  (cf. \cite[ (7)]{ASJ})
\begin{equation}p_{n+m}(x)=\label{TTRRext} S_m^{(n+m)}(x)p_{n}(x)-\l_{n+1}S_{m-1}^{(n+m)}(x)p_{n-1}(x),\end{equation}
which will be used in our proof. 

\begin{theorem}\label{main2}Let $k\in\mathbb{N}_{0}$  and $m\in\{2,3,\ldots,n\}$ be fixed and let $\{p_n\}_{n=0}^{\infty}$ be a sequence of monic polynomials orthogonal on the  interval $(a,b)$ with respect to the weight $w(x)>0$, where the degree of $p_n$ is $n$. Let $c_{2k}$ be an even polynomial of exact degree $2k$. Then the sequence of polynomials $\{g_{n,k}\}_{n=0}^{\infty}$, with degree exactly $n$, 
orthogonal with respect to $c_{2k}(x)w(x)>0$ on $(a,b)$, satisfies
\begin{equation}\label{general44}
c_{2k}(x) g_{n-m,k}(x)=a_{2k-m}(x) p_{n}(x)-G(x) p_{n-1}(x), ~n\in\{2,3,\dots\},\end{equation}
where $a_{2k-m}$ and $G$ are polynomials with $\deg(G)=\max \{m-1,2k-m-1\}$ and 
\begin{equation}\nonumber\deg(a_{2k-m})=\begin{cases}m-2~~\mbox{if}~k\in\{0,1,\ldots,m-1\},\\2k-m~~\mbox{if}~~k\in\{m,m+1,\ldots\}.\end{cases}\end{equation}
 Furthermore, $\deg(G)=m-1$  if and only if 
 $k\in\{0,1,2,\dots,m\}.$
\end{theorem}
\begin{proof}
Fix $k\in\mathbb{N}_{0}$  and $m\in\{2,3,\ldots,n\}$. We apply Christoffel's formula \cite[Thm 2.7.1]{Ism}, see also  \cite[\S 2.5]{szego}, to $g_{n-m,k}$, to obtain
\begin{align}\nonumber
U_{n-m,2k}&c_{2k}(x)g_{n-m,k}(x)\\&=\left|\begin{matrix}\label{det2}
                   p_{n-m}(x_1)& p_{n-m+1}(x_1) & \dots & p_{n-m+2k}(x_1)\\
                   ~~~ p_{n-m}(-x_1)& ~~~p_{n-m+1}(-x_1) & \dots & ~~~p_{n-m+2k}(-x_1)\\
                    p_{n-m}(x_2)& p_{n-m+1}(x_2) & \dots & p_{n-m+2k}(x_2)\\
                     ~~~p_{n-m}(-x_2)& ~~~p_{n-m+1}(-x_2) & \dots & ~~~p_{n-m+2k}(-x_2)\\
                    \dots & \dots & \dots & \dots  \\
                    p_{n-m}(x_k) &  p_{n-m+1}(x_k) & \dots & p_{n-m+2k}(x_k)\\
                    ~~~p_{n-m}(-x_k) &  ~~~p_{n-m+1}(-x_k) & \dots & ~~~p_{n-m+2k}(-x_k)\\
                     p_{n-m}(x) & p_{n-m+1}(x) & \dots & p_{n-m+2k}(x)
                    \end{matrix}\right|\\
               &= \label{det22} \sum_{j=0}^{2k}(-1)^{j} U_{n-m,j} p_{n-m+j}(x),
               \end{align}
               
where $\pm x_i,i\in\{1,2,\dots,k\}$ are the zeros of the even function $c_{2k}(x)$ and $U_{n-m,j}$,  $j\in\{0,1,\dots,2k\}$ is the $2k\times 2k$ determinant obtained from the  $(2k+1)\times (2k+1)$ matrix in (\ref{det2}), by removing the $(2k+1)^{th}$ row and the $(j+1)^{th}$ column, for $j\in\{0,1,\dots,2k\}$. 

Next, we write $U_{n-m,2k}c_{2k}(x)g_{n-m,k}(x)$, and therefore $ p_{n-m+j}(x)$, in terms of $p_n$ and $p_{n-1}$. When $j=m-1$ and $j=m$ we have  $p_{n-m+j}=p_{n-1}$ and $p_{n-m+j}=p_n$ respectively, so we consider the cases where  $j\in\{0,1,\dots,m-2\}$ and $j\in\{m+1,m+2,\dots,2k\}$ separately.  We use (\ref{Beardon}) to express the $(m-1)$ polynomials $p_{n-m+j}$, $j\in\{0,1,\dots,m-2\},$ in terms of  $p_n$ and $p_{n-1}$, i.e., we have, for $j\in\{0,1,\dots,m-2\}$,
\begin{equation}\label{E1} p_{n-m+j}(x)=\frac{S_{m-j-1}^{(n)}(x)}{\l_{n}\ldots \l_{n-(m-j-2)}}p_{n-1}(x)-\frac{S_{m-j-2}^{(n-1)}(x)}{\l_{n}\ldots \l_{n-(m-j-2)}}p_{n}(x).\end{equation} 
Replacing $m$ by $-m+j$ in (\ref{TTRRext}), yields an expression for
$p_{n-m+j}$ when \newline $j\in\{m+1,m+2,\dots,2k\},$ namely
\begin{equation}p_{n-m+j}(x)=\label{E2} S_{-m+j}^{(n-m+j)}(x)p_{n}(x)-\l_{n+1}S_{-m+j-1}^{(n-m+j)}(x)p_{n-1}(x).\end{equation}

Collecting the coefficients of  $p_n$ and $p_{n-1}$ from (\ref{E1}) and (\ref{E2}), and substituting them  into (\ref{det22}) for the different values of $j$, we obtain 
$$U_{n-m,2k}c_{2k}(x)g_{n-m,k}(x)= R^*(x)p_{n}(x) +  G^*(x)p_{n-1}(x),$$
i.e., $$c_{2k}(x)g_{n-m,k}(x)= R(x)p_{n}(x) +  G(x)p_{n-1}(x),$$
with
\begin{eqnarray*}
R(x)&=&\sum_{j=0}^{m-2}\frac{(-1)d_j~S_{m-j-2}^{(n-1)}(x)}{\lambda_n \lambda_{n-1}\ldots \lambda_{n-(m-j-2)}}+\sum_{j=m}^{2k}d_j~S_{j-m}^{(n-m+j)}(x)~=~\frac{R^*(x)}{U_{n-m,2k}},\\
G(x)&=&\sum_{j=0}^{m-2}\frac{d_j~S_{m-j-1}^{(n)}(x)}{\lambda_n\lambda_{n-1}\ldots\lambda_{n-(m-j-2)}}
+d_{m-1}-\lambda_{n+1}\sum_{j=m+1}^{2k}d_j~S_{j-m-1}^{(n-m+j)}(x)\\&=&\frac{G^*(x)}{U_{n-m,2k}}
\end{eqnarray*}
and $d_j=\displaystyle\frac{(-1)^{j}U_{n-m,j}}{U_{n-m,2k}},j\in\{0,1,\dots,2k\}$.  Since $d_1 \neq 0$ (cf. \cite[Thm 2.7.1]{Ism}) and $d_{2k}=1$, it is clear that $\deg(G)=\max \{m-1,2k-m-1\}$. Furthermore, $\deg(R)=\max\{m-2,2k-m\}$ and if
we let $a_{2k-m}(x)=R(x)$, we have \begin{equation}\nonumber\deg(a_{2k-m})=\begin{cases}m-2~~\mbox{if}~k\in\{0,1,\ldots,m-1\},\\2k-m~~\mbox{if}~~k\in\{m,m+1,\ldots\}.\end{cases}\end{equation}
Lastly, since $\deg(G)=\max \{m-1,2k-m-1\}$, it follows directly that $\deg(G)=m-1$ if and only if $k\in\{0,1,2,\dots,m\}.$
\end{proof}

\begin{remark}
\begin{itemize}
    \item[]
    \item[(i)] 
    Equation \eqref{Beardon} corrects \cite[(10)]{Beardon}. In the notation of \cite{Beardon}, equation (10) in \cite[Theorem 4]{Beardon} should read as
$$S_{m-1}(x)p_{n+m}(x)=S_{m}(x)p_{n+m-1}(x)+\l_n \l_{n+1}\dots \l_{n+m-2}p_n(x).$$ The proof of the theorem is correct except for a sign error in the determinants. From \cite[(11)]{Beardon}, it follows that $$S_2(x)V_2(x)-U_2(x)T_2(x)=\l_n.$$
                Substituting $m=2$ in \cite[(12)]{Beardon} and taking determinants yields
     $$S_3(x)V_3(x)-U_3(x)T_3(x)=\l_n\l_{n+1}.$$
                Using induction as discussed in \cite{Beardon}, 
               $$\left|\begin{matrix}
                   S_m(x)& T_m(x)\\
                  U_m(x) & V_m(x)
                \end{matrix}
                \right|=S_m(x)V_m(x)-U_m(x)T_m(x)=\l_n\l_{n+1}\dots \l_{n+m-2},$$ and this corrects the misprint in the proof of \cite[Thm 4]{Beardon}.  

\item[(ii)]
Equation \eqref{general44} in Theorem \ref{main2}  differs from the mixed recurrence equation in \cite[Thm 3]{ASJ}, where the polynomial $g_{n-m,k}$, obtained from making a parameter shift of $k$ units, denotes a polynomial orthogonal with respect to the weight $c_k(x)w(x)$ on $(a,b)$ and $c_k(x)>0$ is a polynomial of degree $k$ in $x$. In \cite[Thm 3]{ASJ}, the parameter shift is equal to the degree of the perturbation whereas Theorem \ref{main2} pertains to the case where the degree of the perturbation is double the parameter shift. 
\end{itemize}
\end{remark}
\section{Meixner-Pollaczek polynomials}
Monic Meixner-Pollaczek polynomials are defined by \cite[\S 9.7]{KLS}
\begin{equation*}
P_n^{(\lambda)}(x;\phi)\label{MPdef}=  i^n(2\lambda)_n\left(\frac{ e^{2i\phi}}{e^{2i\phi}-1}\right)^n
\hypergeom{2}{1}{-n,\lambda+ix}{2\lambda}{1-\frac{1}{e^{2i\phi}}}
\;
\end{equation*}
and are orthogonal with respect to the continuous weight \begin{equation*}\label{MPweight}w(x)=|\Gamma(\lambda+ix)|^2\,e^{(2\phi-\pi)x}\end{equation*} on the interval $(-\infty,\infty),$ for $n\in\nn,$  $\lambda>0$ and $0<\phi<\pi.$ Here,
\[
\Gamma(z):=\int_0^\infty t^{z-1}e^{-t}dt
\]
denotes the Gamma function \cite[(5.2.1)]{DLMF},
\[
\hypergeom{2}{1}
{\alpha_{1},\alpha_{2}}{\beta}{x}
=
\sum_{k=0}^\infty \frac
{(\alpha_{1})_{k}(\alpha_{2})_{k}}
{(\beta)_{k}}\,\frac{x^k}{k!}
\]
the Gauss hypergeometric series \cite[ (15.2.1)]{DLMF} and
\begin{eqnarray*}(\a)_n& =&(\a)(\a+1)\cdots(\a+n-1)~\mbox{for}~
n\geq1\\
(\a)_0&=&1 ~\mbox{when}~ \a\neq0\end{eqnarray*} the Pochhammer symbol.

\medskip
Meixner-Pollaczek polynomials are closely related to (monic) Meixner polynomials \cite[\S 9.10]{KLS}:
\[
M_n(x;\beta,c)=(\beta)_n\left(\frac{c}{c-1}\right)^n
\hypergeom{2}{1}{-n,-x}{\beta}{1-\frac 1c},
\;
\]
orthogonal with respect to the discrete weight $\rho(x)=\frac{c^x(\beta)_x}{x!}$ on $(0,\infty),$ for $0<c<1$ and $\beta>0$. 

\medskip 
We note that
\begin{equation}P_{n}^{(\lambda)}(x;\phi)=
\label{MTOMP}i^n M_{n}(-\lambda -ix;2 \lambda,e^{2i\phi}),\end{equation}
but that Meixner polynomials are orthogonal with respect to a discrete variable while Meixner-Pollaczek polynomials are orthogonal with respect to a continuous variable.

\medskip 
Monic Meixner polynomials satisfy the three-term recurrence equation \cite[\S 9.10]{KLS}
\begin{eqnarray*}\lefteqn{M_{n}(x;\b,c)}\\&=\left(x+\frac{c (\beta +n-1)+n-1}{c-1}\right)M_{n-1}(x;\b,c)-\frac{c (n-1) (\beta +n-2)}{(c-1)^2 } M_{n-2}(x;\b,c),\end{eqnarray*}
with $M_0(x;\b,c)=1$ and $M_{-1}(x;\b,c)=0$ and, when we substitute $x$ with $-\lambda-ix$, $\b$ with $2\lambda$ and $c$ with $e^{2i\phi}$, multiply by $i^n$ and apply (\ref{MTOMP}), we obtain the three-term recurrence equation for the monic Meixner-Pollaczek polynomials
\begin{equation}P_{n}^{(\lambda)}(x;\phi)=\label {TTRRMP}\left(x+\frac{\lambda +n-1}{\tan\phi}\right) P_{n-1}^{(\lambda)}(x;\phi)-\frac{(n-1)(2\lambda +n-2)}{4\sin^2 \phi} P_{n-2}^{(\lambda)}(x;\phi),\end{equation}
with $P_{1}^{(\lambda)}(x;\phi)=x+\l \cot \phi,$ $P_{0}^{(\lambda)}(x;\phi)=1$ and $P_{-1}^{(\lambda)}(x;\phi)=0.$

\subsection{Connection formulae for Meixner-Pollaczek polynomials}\label{MP1}
In order to find new inner bounds for the extreme zeros of Meixner-Pollaczek polynomials, we require mixed three-term recurrence equations involving polynomials $P_{n}^{(\lambda)}(x;\phi)$, $P_{n-1}^{(\lambda)}(x;\phi)$ and $P_{n-m}^{(\lambda+k)}(x;\phi)$. Let $\lambda>0,~ 0<\phi<\pi,~k,n\in \mathbb{N}_{0}$. Since the parameter-shifted polynomial $P_{n-m}^{(\lambda+k)}(x;\phi)$ is orthogonal with respect to the weight function 
\begin{eqnarray*}&&e^{(2\phi-\pi)x} |\Gamma(\lambda+k+ix)|^2\\&=&|(\lambda+ix)_k|^2 w(x)\\
&=& |\lambda+ix|^2|(\lambda+1)+ix|^2\dots |(\lambda+k-1)+ix|^2 w(x)\\
&=&(\l^2+x^2)\left((\l+1)^2+x^2\right)\dots\left((\l+k-1)^2+x^2\right) w(x)\\
&=&c_{2k}(x)w(x)\end{eqnarray*}
where the factor $c_{2k}(x)$ is even, we know with reference to Theorem \ref{main2}, that $\deg(G)=m-1$ is equivalent to $k\in\{0,1,2, \dots,m\}$. Hence, for $m=2$, there are exactly three equations of type \eqref{general44} that will involve a linear coefficient for $P_{n-1}^{(\lambda)}(x;\phi)$. For the case when $k=0$, there is no parameter shift and the relevant equation is the three-term recurrence equation \eqref{TTRRMP}. For the case $k=1$, $c_{2}(x)=\lambda^2+x^2$ and we will use Christoffel's theorem \cite[Thm 2.7.1]{Ism} to determine the mixed three-term recurrence equation necessary to find inner bounds for the extreme zeros of the polynomial $P_{n}^{(\lambda)}(x;\phi)$
\begin{align*}
\left|\begin{matrix}
P_{n-2}^{(\lambda)}(i\lambda;\phi)&P_{n-1}^{(\lambda)}(i\lambda;\phi)\\
P_{n-2}^{(\lambda)}(-i\lambda;\phi)&P_{n-1}^{(\lambda)}(-i\lambda;\phi)\end{matrix}\right|&(\lambda^2+x^2)P_{n-2}^{(\lambda+1)}(x;\phi)\\
=&\left|\begin{matrix}
                  P_{n-2}^{(\lambda)}(i\lambda;\phi) & P_{n-1}^{(\lambda)}(i\lambda;\phi) & P_{n}^{(\lambda)}(i\lambda;\phi)\\
                  P_{n-2}^{(\lambda)}(-i\lambda;\phi) & P_{n-1}^{(\lambda)}(-i\lambda;\phi) & P_{n}^{(\lambda)}(-i\lambda;\phi)\\
                   P_{n-2}^{(\lambda)}(x;\phi) & P_{n-1}^{(\lambda)}(x;\phi) & P_{n}^{(\lambda)}(x;\phi)\end{matrix}\right|.
                   \end{align*}
        By expanding the determinant on the right hand side, and using the three-term recurrence equation \eqref{TTRRMP} to express $P_{n-2}^{(\lambda)}(x;\phi)$ in terms of $P_n^{(\lambda)}(x;\phi)$ and $P_{n-1}^{(\lambda)}(x;\phi)$, we obtain, for $k=1$,
\begin{align}
\frac{(n-1)\left(\lambda ^2+x^2\right) }{2\lambda +n-1}& P_{n-2}^{(\lambda+1)}(x;\phi)\label{eq1}\\\nonumber &=\left(x+\frac{\l}{\tan \phi}\right)P_{n-1}^{(\lambda)}(x;\phi)-\frac{2 \lambda }{2 \lambda +n-1}P_{n}^{(\lambda)}(x;\phi)\end{align}
It is straightforward to verify this equation by comparing coefficients of $x^n$. Using the same procedure, we obtain, for $k=2$,
\begin{align}&\frac{2 (n-1)\sin ^2\phi  }{(2 \lambda +1)(2 \lambda +n-1)(\lambda +n)}\left(\lambda ^2+x^2\right) \left((\lambda +1)^2+x^2\right) P_{n-2}^{(\lambda+2)}(x;\phi)\label{eq2}\\
&\quad\quad=\nonumber\left(x+\frac{ \lambda  (\lambda +1) }{(\lambda +n)\tan\phi} \right) P_{n-1}^{(\lambda)}(x;\phi)-D(x)P_{n}^{(\lambda)}(x;\phi)\end{align} where
{\small{\[
D(x)=\nonumber\frac{(\lambda)_2 (4 \lambda -(n-1) \cos2 \phi +n+1)}{(2 \lambda +1)(2 \lambda +n-1)(\lambda +n)}+\frac{2(n-1) \left((\lambda +\tfrac12) \sin2 \phi- x \sin ^2\phi\right)}{(2 \lambda +1)(2 \lambda +n-1)(\lambda +n)}x.\]}}
Next we consider the case where $m=3$ for Meixner-Pollaczek polynomials in Theorem \ref{main2} and note that equations connecting $P_{n-1}^{(\lambda)}(x;\phi)$, $P_{n}^{(\lambda)}(x;\phi)$ and $P_{n-3}^{(\lambda+k)}(x;\phi)$, with a quadratic coefficient for $P_{n-1}^{(\lambda)}(x;\phi)$, exist if and only if $k\in\{0,1,2,3\}$. The four connection formulae are derived as before and are given in Appendix \ref{appendix}.

\subsection{Inner bounds for the extreme zeros of Meixner-Pollaczek polynomials}\label{innerbounds}
\subsubsection {Using connection formulae for polynomials with a degree difference of two}
Using  \cite[Cor. 2.2]{DJ}, the zeros of the coefficients of $P^{(\lambda)}_{n-1}(x;\phi)$ obtained from equations (\ref{TTRRMP}), (\ref{eq1}) and (\ref{eq2}) are inner bounds for the extreme zeros of the Meixner-Pollaczek polynomial $P_n^{(\lambda)}(x;\phi)$. Hence, the points \begin{equation}\label{bnds}B_n(k)=-\frac{(\l)_k~(\l)_n}{(\l)_{n+k-1}\tan\phi}, k \in\{0,1,2\},\end{equation} satisfy $x_{1,n}<B(k)<x_{n,n}$ for $k \in\{0,1,2\}$ where $x_{i,n}, i\in\{1,2,\dots,n\}$ denote the zeros of $P_n^{(\l)}(x;\phi).$ We note that for $n=1$, $B_1(k)=-\l \cot \phi=x_{1,1}$ for each $k\in\{0,1,2\}$. It can easily be shown that, for $n\in\{2,3,\dots\}$,
$B_{n}(0)<B_{n}(1)<B_{n}(2)$ for all values of $\lambda>0$ and $0<\phi<\frac{\pi}{2}$ while $B_{n}(2)<B_{n}(1)<B_{n}(0)$ when $\lambda>0$ and $\frac{\pi}{2}<\phi<\pi$. Hence we have, for $\la>0$ and $0<\phi<\frac{\pi}{2}$, that 
\begin{equation}\label{B} x_{1,n}<B_{n}(0)<B_{n}(1)<B_{n}(2)<x_{n,n}\end{equation}
and this is illustrated for various values of $n$, $\la$ and $\phi$ in Table \ref{MeixnerP4}.
\begin{table}[!ht] \caption{Comparison of bounds given in (\ref{bnds}), for the extreme zeros of $P_{n}^{(\lambda)}(x;\phi)$ for different values of $n, \lambda$ and $\phi\in(0,\frac{\pi}{2}$).}\label{MeixnerP4}
\begin{center}
\begin{tabular}{c|c|c|c|c|c|c}
$n$&$\lambda$&$\phi$&$x_{1,n}$&$B_{n}(0)$ & $B_{n}(2)$  &$x_{n,n}$\\
\hline
6&0.89&0.25&-34.28&-23.07&-0.96&-0.63\\
6&0.89&0.1&-86.31&-58.70&-2.43&-2.17\\
6&0.01&0.1&-71.38&-49.93&-0.017&-0.015\\
10&0.5&0.08&-186.87&-118.50&-0.891&-0.768\\
20&0.01&0.1&-323.1&-189.5&-0.005&-0.043\\
20&1.5&0.1&-350.7&-201.8&-1.280&-1.082\\
\end{tabular}
\end{center}
\end{table}
As can be seen from Table \ref{MeixnerP4}, the inner bound $B_n(2)$ is a relatively good lower bound for the largest zero $x_{n,n}$ of $P_n^{(\lambda)}(x;\phi)$ but $B_n(0)$ is not a sharp upper bound for the smallest zero $x_{1,n}$.
When $\l>0$ and $\frac{\pi}{2}<\phi<\pi$, we have that 
\begin{equation}\nonumber x_{1,n}<B_{n}(2)<B_{n}(1)<B_{n}(0)<x_{n,n}\end{equation}
and this is illustrated for various values of $n$, $\la$ and $\phi$ in Table \ref{MeixnerP44444444}.
\begin{table}[!ht] \caption{Comparison of bounds given in (\ref{bnds}), for the extreme zeros of $P_{6}^{(\lambda)}(x;\phi)$ for different values of $\lambda$ and $\phi\in(\frac{\pi}{2},\pi$).}\label{MeixnerP44444444}
\begin{center}
\begin{tabular}{c|c|c|c|c|c}
$\lambda$&$\phi$&$x_{1,6}$&$B_{6}(2)$ & $B_{6}(0)$  &$x_{6,6}$\\
\hline
0.89&1.7&-2.194&0.0037&1.187&4.0569\\
0.89&3.14&142.865&153.289&3698.23&5426.31\\
0.01&3.14&1.05511&1.05518&3145.69&4487.87\\
0.08&2.5&-0.356&0.019&6.8&10.6806\\
\end{tabular}
\end{center}
\end{table}

Next we will use the connection formulae involving polynomials with a higher degree difference in order to improve on the bounds obtained in this section. 
\subsubsection{Using connection formulae for polynomials with a degree difference of three}
As before, by using \cite[Cor. 2.2]{DJ}, we know that the zeros of the coefficients of $P_{n-1}^{(\lambda)}(x;\phi)$ in equations \eqref{rr1}, \eqref{rr2}, \eqref{rr3} and \eqref{rr4}, given in the Appendix, are inner bounds for the extreme zeros of $P_{n}^{(\lambda)}(x;\phi),n\in\{2,3,\dots\}$. We will denote these zeros, which all are inner bounds for the extreme zeros of $P_{n}^{(\lambda)}(x;\phi)$, by $F_n^{\pm}(k)$, $k\in\{0,1,2,3\}$.

The zeros of $E_2(x)$ in \eqref{rr1}, $H_2(x)$ in \eqref{rr2}, $M_2(x)$ in \eqref{rr3} and  $R_2(x)$ in \eqref{rr4} are, respectively,  
\begin{align}\label{z1}F_n^{\pm}(0)=&-( \lambda+ n-\tfrac 32) \cot \phi\pm\tfrac{1}{2 \sin \phi} \sqrt{(n-1)  (n+2 \lambda-2)+\cos^2  \phi }\\ \label{z2} F_n^{\pm}(1)=&-\tfrac{1}{2} (2 \lambda +n-1) \cot (\phi )\pm\tfrac{1}{2 \sin \phi}\sqrt{(n-1) \left(2 \lambda +(n-1) \cos ^2\phi\right)}\\
F_n^{\pm}(2)=&-\tfrac{1}{2} (2 \lambda +1) \cot \phi \label{z3}\\
&\qquad \pm\frac{ \sqrt{-2 \lambda  (4 \lambda +3)+(2 \lambda +n) \cos (2 \phi )+8 \lambda  (\lambda +1) n+n}}{2  \sin \phi\sqrt{2 (2 \lambda +n)}}\nonumber\\
\label{z4}
F_n^{\pm}(3)=&\frac{- (4 \lambda  (\lambda +2)+3) (\lambda +n) \sin (2 \phi)\pm  \frac{1}{2}\sqrt{B}}{\sin^2 \phi ~(4 \lambda  (2 \lambda +1)+6 n^2+12 \lambda  n)}
\end{align}  
where 
\begin{align*}
B &=4 (4 \lambda  (\lambda +2)+3)^2 (\lambda +n)^2 \sin ^2(2 \phi)\\
&-8 \lambda  (\lambda +2) \left(2 \lambda  (2 \lambda +1)+3 n^2+6 \lambda  n\right) \sin ^2\phi\\
&\quad\times\left(4 \lambda ^2+\left(4 \lambda ^2+n^2+2 \lambda  (n+3)+2\right) \cos (2 \phi)-n^2-2 \lambda  (n-5)+4\right).\end{align*}

It is clear that each one of the values $F_n^{-}(i)$ and $F_n^{+}(i),i\in\{0,1,2,3\}$ is a lower bound for the largest and an upper bound for the smallest zero of $P_n^{\l}(x),$  i.e.,
\begin{equation*}x_{1,n}<F_n^{-}(i)<F_n^{+}(i)<x_{n,n} ~~\mbox{for each}~~ i\in\{0,1,2,3\}.\end{equation*}
In the following theorem we prove which one of the values $F_n^{-}(i),i\in\{0,1,2,3\},$ is the strongest upper bound for the smallest zero when $\phi \in (0, \frac {\pi}{2})$ and we will also compare this bound with the best upper bound obtained from connection formulae involving polynomials with a degree difference of two.
\begin{theorem}
For all $n\in\{2,3,\dots\},$ $\l>0$ and $\phi \in (0, \frac {\pi}{2}),$
\begin{itemize}
\item[(a)]
$F_n^-(0)\leq F_n^-(1)\leq F_n^-(2)\leq  F_n^-(3).$
\item [(b)] $F_n^-(0)\leq F_n^-(1)\leq B_n(0)\leq B_n(1)\leq B_n(2).$
\end{itemize}
\end{theorem}
\proof
Let  $n\in\{2,3,\dots\},$ $\l>0$ and $\phi \in (0, \frac {\pi}{2}).$
\begin{itemize}
\item[(a)]  We will divide the proof into 3 parts and prove that
 \begin{itemize}
     \item [(i)] $F_n^-(0)\leq F_n^-(1);$
      \item [(ii)] $F_n^-(1)\leq F_n^-(2);$
       \item [(iii)] $F_n^-(2)\leq F_n^-(3).$
       \end{itemize}
\begin{itemize}
     \item [(i)] From  \eqref{z1} 
 and \eqref{z2}, we obtain
\begin{align*}&\frac{F_n^-(1)}{F_n^-(0)}\\&=\frac{(2 \lambda +n-1) \cos (\phi )+\sqrt{2 \lambda  (n-1)+(n-1)^2 \cos ^2(\phi )}}{(2 \lambda +2 n-3) \cos (\phi )+\sqrt{(n-1) (2 \lambda +n-2)+\cos ^2(\phi )}}\\
&=\frac{(2 \lambda +n-1) \cos (\phi )+\sqrt{2 \lambda  (n-1)+(n-1)^2 \cos ^2(\phi )}}{(2 \lambda +n-1) \cos (\phi )+(n-2)\cos{\phi}+\sqrt{2\l(n-1)+(n-2)_2+\cos ^2(\phi )}}\\
&=\frac{(2 \lambda +n-1) \cos (\phi )+\sqrt{2\l(n-1) +(n-1)^2 \cos ^2(\phi )}}{(2 \lambda +n-1) \cos (\phi )+\sqrt{(n-2)^2\cos^2{\phi}}+\sqrt{2\l(n-1) +(n-2)_2+\cos ^2(\phi )}}\\
&\leq\frac{(2 \lambda +n-1) \cos (\phi )+\sqrt{2\l(n-1) +(n-1)^2 \cos ^2(\phi )}}{(2 \lambda +n-1) \cos (\phi )+\sqrt{(n-2)^2\cos^2{\phi}+2\l(n-1) +(n-2)_2+\cos ^2(\phi )}}\\&\qquad \left(\text{since }\sqrt{a+b}\leq\sqrt{a}+\sqrt{b}\implies\frac{1}{\sqrt{a}+\sqrt{b}}\leq \frac{1}{\sqrt{a+b}}\right)\\
&=\frac{(2 \lambda +n-1) \cos (\phi )+\sqrt{2\l(n-1)+(n^2-4n+5) \cos ^2(\phi )+2(n-2)\cos^2\phi}}{(2 \lambda +n-1) \cos (\phi )+\sqrt{2\l(n-1) +(n^2-4n+5)\cos^2{\phi}+(n-1)(n-2)}}\\
&\leq \frac{(2 \lambda +n-1) \cos (\phi )+\sqrt{2\l(n-1) +(n^2-4n+5) \cos ^2(\phi )+2(n-2)}}{(2 \lambda +n-1) \cos (\phi )+\sqrt{2\l(n-1) +(n^2-4n+5)\cos^2{\phi}+(n-1)(n-2)}}\\
&\leq 1,
\end{align*} therefore $F_n^-(1)\ge F_n^-(0)$, since $F_n^-(0)<0.$ 
\item[(ii)] From  \eqref{z2} 
 and \eqref{z3}, we obtain
\begin{align*}\frac{F_n^-(2)}{F_n^-(1)}&=\frac{(2 \lambda +1) \cos (\phi )+\sqrt{\frac{4 \lambda  (\lambda +1) (n-1)}{2 \lambda +n}+\cos ^2(\phi )}}{(2 \lambda +n-1) \cos (\phi )+\sqrt{2 \lambda  (n-1)+(n-1)^2 \cos ^2(\phi )}}\le 1,\end{align*}
 since $$\frac{\frac{4 \lambda  (\lambda +1) (n-1)}{2 \lambda +n}}{2 \lambda  (n-1)}=\dfrac{2\l+2}{2\l+n}\le 1$$ for each $n\in\{2,3,\dots\},$ and, taking into account that $F_n^-(1)$ is negative, we can conclude that $F_n^-(2)\ge F_n^-(1).$ 
\item [(iii)] From  \eqref{z2} 
 and \eqref{z3}, we obtain, after simplifying, that
\begin{align}
&\frac{F_n^-(3)}{F_n^-(2)}= \label{e} \frac{2 B (\lambda +n) \cos \phi}{A(2 \lambda +1) \cos \phi   +\sqrt{ A^2\cos ^2\phi +\frac{A^24 (\lambda +1) \lambda  (n-1)}{(2 \lambda +n)}}}\\
+&\nonumber \frac{\sqrt{4 \cos ^2\phi \left(B^2 (\lambda +n)^2-AC \lambda  (\lambda +2) \right)+ 4 A (\lambda +2) \lambda  (n-1) (2 \lambda +n+1)}}{A(2 \lambda +1) \cos \phi   +\sqrt{ A^2\cos ^2\phi +\frac{A^24 (\lambda +1) \lambda  (n-1)}{(2 \lambda +n)}}}\end{align}
where 
\begin{align*}
A&= \lambda ^2+2 \lambda +3 n^2+6 \lambda  n\\
B&=4 \lambda ^2+8 \lambda +3\\
C&=4 \lambda ^2+6 \lambda +n^2+2 \lambda  n+2.
\end{align*}
For $n\in\{2,3,\dots\}$ and $\l>0$ we have
\begin{equation}\label{e1}0<\frac{2 B (\lambda +n)}{A (2 \lambda +1)}=\frac{2 (2 \lambda +3) (\lambda +n)}{2 \lambda  (2 \lambda +1)+3 n^2+6 \lambda  n}\leq1,\end{equation} since \[2 (2 \lambda +3) (\lambda +n)-\left(2 \lambda  (2 \lambda +1)+3 n^2+6 \lambda  n\right)=-(n-2) (2 \lambda +3 n)\leq0,\] while
\begin{align}\label{e2}0&<\frac{4 A \lambda  (\lambda +2) (n-1) (2 \lambda +n+1)}{\frac{A^2 (4 \lambda  (\lambda +1) (n-1))}{2 \lambda +n}}\\
&=\nonumber\frac{(\lambda +2) (2 \lambda +n) (2 \lambda +n+1)}{A (\lambda +1)}\\
&\leq 1\nonumber,\end{align} since 
\[(\lambda +2) (2 \lambda +n) (2 \lambda +n+1)-A(\lambda +1) =(2 \lambda +1) (2-n) (\lambda +n)\leq 0\] and
\begin{equation}\label{e3}\frac{4 \left(B^2 (\lambda +n)^2-A C \lambda  (\lambda +2)\right)}{A^2}\leq 1,\end{equation} since
\begin{align*}
 4 &\left(B^2 (\lambda +n)^2-A F \lambda  (\lambda +2)\right)-A^2\\
 &=-3(4 \lambda^2  +8\lambda +3) n(n-2)(2 \lambda +n) (2 \lambda +n+2)\\
 &\leq 0.   
\end{align*}
Using \eqref{e}, \eqref{e1}, \eqref{e2} and \eqref{e3} and taking into account the fact that $F_n^-(2)<0$, we see that, for all $n\in\{2,3,\dots\},$ $\l>0$ and $\phi \in (0, \frac {\pi}{2})$,  $F_n^-(3)\geq F_n^-(2)$.
\end{itemize}
\item[(b)] Using the formulas for $B_n(0)$, $F_n^{-}(0)$ and $F_n^{-}(1)$ in \eqref{bnds}, \eqref{z1} and \eqref{z2}, respectively, we obtain
\[B_n(0)-F_n^{-}(0)= \frac{\sqrt{\cos^2 \phi+ (n-1)(n+2\lambda-2)} -\cos \phi}{2\sin \phi}>0\]  and
\[B_n(0)-F_n^{-}(1)=\frac{\sqrt{(n-1)^2\cos^2 \phi+ 2\lambda(n-1)}-(n-1)\cos \phi}{2\sin \phi}>0\] and, using also \eqref{B}, we conclude that, for $n\in \mathbb{N}$, $\lambda>0$ and $\phi\in(0,\tfrac{\pi}{2})$, $$F_n^-(0)\leq F_n^-(1)\leq B_n(0)\leq B_n(1)\leq B_n(2).$$\endproof
\end{itemize}

Hence, for $n\in\{2,3,\dots\},$ $\l>0$ and $\phi \in (0, \frac {\pi}{2})$, $F_n^{-}(0)$ provides the best new upper bound for the smallest zero of  $P_n^{(\l)}(x;\phi)$ and it also is a sharper upper bound than the bound $B_n(0)$ obtained in \eqref{bnds}. Furthermore, $F_2^-(0)=F_2^-(1)=F_2^-(2)=F_2^-(3)$. 

\medskip
We note that for each $n\in\{2,3,\dots\}, \l>0$ and $\phi \in (\frac {\pi}{2}, \pi)$,$$F_n^+(0)\leq F_n^+(1)\leq F_n^+(2)\leq  F_n^+(3).$$ 
We can therefore show that, for $\phi \in (\frac {\pi}{2}, \pi)$, $F_n^{+}(3)$ is the best new lower bound for the largest  zero of  $P_n^{\l}(x;\phi)$ which also is a sharper lower bound than $B_n(2)$ obtained in \eqref{bnds}.

We illustrate the sharpness of the obtained bounds in Table \ref{MeixnerP44}.
\begin{table}[!ht] \caption{
Comparison of bounds given in (\ref{bnds}), \eqref{z1} and \eqref{z4} for the extreme zeros of $P_{n}^{(\lambda)}(x;\phi)$ for different values of $n, \lambda$ and $\phi\in(0,\frac{\pi}{2}$).}\label{MeixnerP44}
\begin{center}
{\small{\begin{tabular}{c|c|c|c|c|c|c|c|c}
$n$&$\lambda$&$\phi$&$x_{1,n}$&$ F_n^{-}(0)$&$B_{n}(0)$ & $B_{n}(2)$ & $F_n^{+}(3)$&$x_{n,n}$\\
\hline
6&0.89&0.25&-34.28&-32.15&-23.07&-0.96&-0.662 &-0.63\\
6&0.89&0.1&-86.31&-81.10&-58.70&-2.43&-2.174&-2.167\\
6&0.01&0.1&-71.38&-67.95&-49.93&-0.0167&-0.01482&-0.01478\\
10&0.5&0.08&-186.87&-168.91&-118.50&-0.891&-0.773&-0.768\\
20&0.01&0.1&-323.1&-277.29&-189.5&-0.005&0.0015&0.0043\\
20&1.5&0.1&-350.69&-299.59&-201.8&-1.280&-1.169&-1.082\\
\end{tabular}}}
\end{center}
\end{table}

Finally, we will consider how closely the new inner and known outer bounds for the extreme zeros surround the extremal zeros as a way to measure the sharpness of the bounds. 

Upper bounds for the largest zeros and lower bounds for the smallest zeros of $P_n^{(\lambda)}(x;\phi)$ were proved by Ismail and Li in \cite[Thm 7]{Ismail_Li_1992}. 
\begin{theorem}[\cite{Ismail_Li_1992}] \label{IL} Let $n \in \nn$ be fixed, $\lambda>0$ and $0<\phi<\pi$. Denote the zeros of $P_n^{(\lambda)}(x;\phi)$ by $x_{1,n}<x_{2,n}<\dots<x_{n,n}$, then
\begin{itemize}
\item[(i)] $x_{1,n}>U=(n+\lambda-\tfrac 32)\cot \phi-\frac{1}{2\sin\phi}\sqrt{\cos^2 \phi+4\ (n-1)(n+2\lambda-2)}$
\item[(ii)] $x_{n,n}<V=(n+\lambda-\tfrac 32)\cot \phi+\frac{1}{2\sin\phi}\sqrt{\cos^2 \phi+4\ (n-1)(n+2\lambda-2)}$
\end{itemize} \end{theorem}
\proof Applying \cite[Thms 2 and 3]{Ismail_Li_1992} to the three-term recurrence equation \eqref{TTRRMP} (with $n$ replaced by $n+1$) by letting $\alpha=1$, $\beta=\frac{\lambda+n}{\tan\phi}$ and $\lambda_n=\frac{n(n+2\lambda-1)}{4\sin^2\phi}$ yields the result. 
\endproof
It is clear from numerical experiments that the interior bounds obtained for extreme zeros of Meixner-Pollaczek polynomials in this paper are much sharper than the exterior bounds obtained by Ismail and Li in \cite{Ismail_Li_1992} as illustrated by the examples in Table \ref{MeixnerP444}.
\begin{table}[!ht] \caption{Comparison of inner bounds given in (\ref{z1}) and \eqref{z4} with outer bounds given in Theorem \ref{IL} for the extreme zeros of $P_{n}^{(\lambda)}(x;\phi)$ for different values of $n, \lambda$ and $\phi\in(0,\frac{\pi}{2}$).}\label{MeixnerP444}
\begin{center}
\begin{tabular}{c|c|c|c|c|c|c|c|c}
$n$&$\lambda$&$\phi$&$U$&$x_{1,n}$&$F_n^{-}(0)$&$F_n^{+}(3)$&$x_{n,n}$&$V$\\
\hline
6&0.89&0.25&-214&-34.28&-32.15&-0.662 &-0.63&257\\
6&0.89&0.1&-530&-86.31&-81.10&-2.174&-2.167&638\\
6&0.01&0.1&-363&-71.38&-67.95&-0.01482&-0.01478&453\\
10&0.5&0.08&-1921&-186.87&-168.91&-0.773&-0.768&2146\\
20&0.01&0.1&-6680&-323.1&-277.29&0.0015&0.0043&7048\\
20&1.5&0.1&-7799&-350.69&-299.59&-1.169&-1.082&8198\\
\end{tabular}
\end{center}
\end{table}
\section{Pseudo-Jacobi polynomials}\label{PJ}
For $a,b$ real, the monic Pseudo-Jacobi polynomials \cite{JT,KLS} are defined by 
\begin{equation*}
    P_n(x;a,b)=\frac{2^n (a+ib+1)_n} {i^n (2 a+n+1)_n}~\hypergeom{2}{1}{-n,2 a+n+1}{a+ib +1}{\frac{1-i x}{2}}, \end{equation*}
$n\in\{0,1,2\dots\}$ and, for $b\in\R$ and $a<-n$, they are orthogonal on $\R$ with respect to the weight function $w(x)=(1+x^2)^a e^{2b \arctan x}.$ 

\smallskip For $b\in\R,$ $a<-n$  and $k\in\nn_{0}$ such that $a+k<-n,$ the shifted polynomials \newline$P_n(x;a+k,b)$ are orthogonal on the real line with respect to the weight function
$$(1+x^2)^{k+a} e^{2b \arctan x}=(1+x^2)^{k} w(x)=c_{2k}(x)w(x)>0.$$ 
By Theorem \ref{main2}, Pseudo-Jacobi polynomials satisfy mixed three-term recurrence equations of the form (\ref{general44}), involving $P_{n}(x;a,b),P_{n-1}(x;a,b)$ and $P_{n-2}(x;a+k,b)$, where the coefficient of $P_{n-1}(x;a,b)$ is linear, if and only if $k\in\{0,1,2\}.$

Applying Theorem \ref{main2} to the three-term recurrence equation satisfied by monic Pseudo-Jacobi polynomials (cf. \cite[(3)]{JT})
\begin{align*}
P_{n}(x;a,b)=&\left(x+\frac{a b}{(a+n-1) (a+n)}\right) P_{n-1}(x;a,b)\\&+\frac{(n-1) (2 a+n-1) \left((a+n-1)^2+b^2\right)}{(a+n-1)^2 \left(4 (a+n-1)^2-1\right)}P_{n-2}(x;a,b),
\end{align*}
the equation involving a parameter shift of one unit,
\begin{eqnarray*}
 \label{PSJplus1}P_{n}(x;a,b)&=&\frac{2 a+n}{2 a+2 n-1}\left(x+\frac{b}{a+n}\right) P_{n-1}(x;a,b)\\
 &+&\frac{(n-1) \left(x^2+1\right) }{2 a+2 n-1} P_{n-2}(x;a+1,b)  
\end{eqnarray*}
and the mixed three-term recurrence equation that involves the largest possible parameter shift, namely $k=2$, we obtain (cf. \cite[Thm 2.3]{JT})
\begin{equation}\label{bounds}
 D_n(0)=\frac{-ab}{(a+n)(a+n-1)},D_n(1)=\frac{-b}{a+n}\mbox{~and~}D_n(2)=\frac{-b}{a+1},   \end{equation}
respectively, such that \begin{equation}\label{BPSJ2}y_{1,n}<D_n(i)<y_{n,n}, i=0,1,2,\end{equation} where $y_{i,n}, i\in\{1,2,\dots,n\}$ denote the zeros of $P_n(x;a,b)$ in ascending order. We note that $$P_1(x;a,b)=x+\frac{b}{a+1}$$
and hence $D_n(2)=y_{1,1}$ for each $n\ge 1$.

Furthermore, for $n\geq 2,$  $$y_{1,n}<D_n(2)\leq D_n(1)\leq D_n(0)<y_{n,n}$$ if $b\geq0$, $a<-n$. The order of these bounds will reverse when $b\leq 0$ and $a<-n$. In Tables \ref{PSJbounds2} and \ref{PSJbounds1} we illustrate the sharpness of these bounds when $b\geq 0$. 

\begin{table}[!ht] \caption{\small{Comparison of bounds (given in (\ref{bounds})) for the extreme zeros of $P_{n}(x;a,b),b\geq 0,$ for different values of $n, a$ and $b$.}}\label{PSJbounds2}
\begin{center}
\begin{tabular}{c|c|c|c|c|c|c}
$n$&$a$&$b$&$y_{1,n}$&$D_{n}(2)$ 
& $D_{n}(0)$ &$y_{n,n}$\\
\hline
8&-8.0016&20&1.78751&2.85649&99860.2&99860.3\\
7&-7.01&50&5.62006&8.31947&34703&34704.1\\
16&-17.5&13&0.008285&0.787879&60.6667&69.123\\
\end{tabular}
\end{center}
\end{table}

In the specific case when $b=0$, the weight function of Pseudo-Jacobi polynomials reduces to an even function of $x$ and the zeros of the polynomial $P_n(x;a,0)$ are symmetric about the origin, which means that all polynomials of odd degree will have $x=0$ as a zero.  This case, where $b=0$ is illustrated in Table \ref{PSJbounds1}. 

\begin{table}[!ht] \caption{\small{Comparison of bounds (in (\ref{bounds})) for the extreme zeros of $P_5(x;a,b),b\ge0,$ for different values of $a$ and $b$.}}\label{PSJbounds1}
\begin{center}
{\small\begin{tabular}{c|c|c|c|c|c|c}
$a$&$b$&$y_{1,5}$&$D_5(2)$ 
 & $D_5(1)$  & $D_5(0)$ & $y_{5,5}$\\\hline
-10&8&0.3455&0.8889&1.6&2.6667&3.5733\\
-5.5&8&1.1189&1.7778&16&58.6667&60.7767\\
-5.0015&3&0.245139&0.749719&2000&9988.02&9988.03\\
-5.5&0&-2.1428&0&0&0&2.1428\\
\end{tabular}}
\end{center}
\end{table}

In Table \ref{PSJbounds3} we also show the upper  bound for the largest zero of Pseudo-Jacobi polynomials, as given in \cite[Thm 2.4]{JT}.  
{\small{\begin{table}[!ht] \caption{\small{Comparison of bounds (given in (\ref{bounds}) and in \cite[Thm 2.4 (i)]{JT}) for the extreme zeros of $P_{n}(x;a,b),b\geq 0,$ for different values of $n, a$ and $b$.}}\label{PSJbounds3}
\begin{center}
{\small{\begin{tabular}{c|c|c|c|c|c|c|c}
$n$&$a$&$b$&$y_{1,n}$&$D_{n}(2)$ 
& $D_{n}(0)$ &$y_{n,n}$&$x_n$ in \cite{JT}\\
\hline
8&-8.0016&20&1.78751&2.85649&99860.2&99860.3&99860.5\\
7&-7.01&50&5.620&8.319&34703&34704.1&34706.8\\
16&-17.5&15&0.138&0.909&70.00&79.680&96.77\\
\end{tabular}}}
\end{center}
\end{table}}}
\begin{remark}
\begin{itemize}
\item[]
\item[(i)] The formula for $D_n(0)$ in \eqref{bounds} corrects the typographical error for the bound $D_n(0)$ in \cite[Thm 2.3]{JT}. 
\item[(ii)] Monic Pseudo-Jacobi polynomials are connected to the Jacobi polynomials  \cite[p. 221]{KLS} and
$$P_n(x;a,b)=\frac{(-2 i)^n n!}{(2 a+n+1)_n} P_n^{(a+b i,a-i b)}(i x).$$
The bounds $iB_n(0)$ and $iB_n(2)$ will therefore be inner bounds for the extreme zeros of the Jacobi polynomial $P_n^{(a+b i,a-i b)}(ix)$. This can also be seen from the interlacing of the zeros of $(x-\frac{\beta-\alpha}{\alpha +
 \beta+2})P_{n-1}^{(\alpha +2,\beta+2)}(x)$ and $P_{n+1}^{(\alpha,\beta)}(x)$, proved in \cite[Thm 2.3]{DJJ}, by using the equation
 \begin{equation}\label{JP}c_4(x)P_{n-2}^{(\alpha+2,\beta+2)}(x)=\left(x - \frac{\beta-\alpha}{\alpha +
 \beta+2}\right)P_{n-1}^{(\alpha,\beta)}(x)+A_2(x)P_{n}^{(\alpha,\beta)}(x),\end{equation} where $c_4(x)=(1-x^2)^2$ and $A_2(x)$ is a polynomial of degree $2$ in $x$. The extra interlacing point $\displaystyle{\frac{\beta-\alpha}{\alpha + \beta+2}}$ is equal to the bound $i B_n(2)$, when we replace $\alpha$ with $a+ib$ and  $\beta$ with $a-ib$.  
The mixed three-term recurrence equation (\ref{JP}) is, however, not of the form of (\ref{general44}) and we can apply \cite[Thm 3]{ASJ} in this case. By shifting $\alpha$ and $\beta$ optimally (but separately), namely by $4$ units, will provide us with more precise bounds,  cf. \cite{DJ}. \end{itemize}
\end{remark}
\section*{Acknowledgements}
The authors gratefully acknowledge support by the National Institute for Theoretical and Computational Sciences (NITheCS) towards a writing retreat. They also thank the anonymous referees for useful comments that helped to improve the manuscript.

\section*{Appendix}\label{appendix}
For $n\ge2,$ the recurrence equation involving $P_{n}^{(\lambda)}(x;\phi),$ $P_{n-1}^{(\lambda)}(x;\phi)$ and $P_{n-3}^{(\lambda)}(x;\phi),$ is
\begin{align}\label{rr1}
&\frac{(n-2)_2  (2 \lambda +n-3)_2}{4 \sin ^2\phi} P_{n-3}^{(\lambda)}(x;\phi)\\
&=E_2(x)P_{n-1}^{(\lambda)}(x;\phi)-\nonumber4 \sin \phi  ((\lambda +n-2) \cos \phi+x \sin \phi )P_{n}^{(\lambda)}(x;\phi),\end{align}
where
\begin{align*}E_2(x)&=n^2+(2 \lambda -3) n+2 \left((\lambda -1)^2+x^2\right)\\
+&2 \cos (2 \phi ) \left(\lambda ^2-3 \lambda +2 \lambda  n+(n-3) n-x^2+2\right)+2 x (2 \lambda +2 n-3) \sin (2 \phi )\end{align*} 
The recurrence equation involving $P_{n}^{(\lambda)}(x;\phi),$ $P_{n-1}^{(\lambda)}(x;\phi)$ and $P_{n-3}^{(\lambda+1)}(x;\phi),$ is
\begin{align}\label{rr2}\frac{1}{2} (n-2)_2 &\left(\lambda ^2+x^2\right) P_{n-3}^{(\lambda+1)}(x;\phi)\\\nonumber&=H_2(x)P_{n-1}^{(\lambda)}(x;\phi)-2 \sin \phi (\lambda  \cos \phi+x \sin \phi ) P_{n}^{(\lambda)}(x;\phi),\end{align}
where
\begin{equation}
H_2(x)=\lambda ^2+\cos (2 \phi) \left(\lambda  (\lambda +n-1)-x^2\right)+x (2 \lambda +n-1) \sin (2 \phi)+x^2.
\end{equation}
The recurrence equation involving $P_{n}^{(\lambda)}(x;\phi),$ $P_{n-1}^{(\lambda)}(x;\phi)$ and $P_{n-3}^{(\lambda+2)}(x;\phi),$ is
\begin{align}\label{rr3}2 (n-2)_2 &\left(\lambda ^2+x^2\right) \left((\lambda +1)^2+x^2\right) \sin \phi   P_{n-3}^{(\lambda+2)}(x;\phi)\\
&=\nonumber M_2(x)\frac{2 \lambda +n-1}{\sin \phi}  P_{n-1}^{(\lambda)}(x;\phi)\\
&-\nonumber 4 (2 \lambda +1) ((\lambda)_2   \cos \phi+x (\lambda +n-1) \sin \phi )P_{n}^{(\lambda)}(x;\phi),\end{align}
where 
\begin{align*}M_2(x)&=(2 \lambda +n)\left(\left(\lambda ^2+\lambda -x^2\right) \cos (2 \phi )+(2 \lambda +1) x \sin (2 \phi )\right)\\
&+n \left(x^2-\lambda  (\lambda +1)\right)+2 \lambda  \left((\lambda +1)^2+x^2\right)\end{align*}
 Finally, the recurrence equation involving $P_{n}^{(\lambda)}(x;\phi),$ $P_{n-1}^{(\lambda)}(x;\phi)$ and $P_{n-3}^{(\lambda+3)}(x;\phi),$ is
\begin{align}\label{rr4}&4 (n-2)_2 \left(\lambda ^2+x^2\right) \left((\lambda +1)^2+x^2\right) \left((\lambda +2)^2+x^2\right) \sin ^3\phi  P_{n-3}^{(\lambda+3)}(x;\phi)\\
&=\nonumber \frac{(\lambda +1) (2 \lambda +n-1)}{\sin \phi} R_2(x) P_{n-1}^{(\lambda)}(x;\phi)-2G_2(x) P_{n}^{(\lambda)}(x;\phi),\end{align}
where 
{\small{\begin{align*}R_2(x)&=-n^2 \left(\lambda  (\lambda +2)-3 x^2\right)+n \left(\lambda  (\lambda +2)-3 x^2\right) ((2 \lambda +n) \cos (2 \phi )-2 \lambda )\\
&+2 (2 \lambda +1) (2 \lambda +3) x (\lambda +n) \sin (2 \phi )\\
&+2 \lambda  (2 \lambda +1) \left((\lambda +1) (\lambda +2)-x^2\right) \cos (2 \phi)+2 \lambda  (2 \lambda +1) \left((\lambda +2)^2+x^2\right)\end{align*}}}
and 
{\small{\begin{align*}
G_2(x)&=x \sin \phi(n^2-(n-2)_2\cos (2 \phi ))\left(3 \lambda  (\lambda +2)-x^2+2\right)\\
&+ n x \sin \phi \left(16\lambda^3+39\lambda^2+26\lambda+ 3 x^2+6-2 \left(\lambda  (13-2 \lambda  (2 \lambda  (\lambda +1)-5))+x^2+4\right)\right)\\
&+(\lambda +1) \cos \phi \\
&\quad\times\left(\lambda  (\lambda +2) (8 \lambda  (\lambda +2)+(n-2)_2  \cos (2 \phi)-n^2+3n+4)+6 (n-2)_2  x^2 \sin ^2\phi\right).
\end{align*}}}

\begin{thebibliography}{99}

\bibitem{Area_et_al_2004} I. Area, D.K. Dimitrov, E. Godoy and A. Ronveaux, Zeros of Gegenbauer and Hermite Polynomials and Connection Coefficients, {\it{Math. Comp.}} 73 (2004), 1937--1951.

\bibitem{Area_et_al_2012} I. Area, D.K. Dimitrov, E. Godoy and F.R. Rafaeli, Inequalities for zeros of Jacobi polynomials via Obrechkoff's theorem, {\it{Math. Comp.}} 81 (2012), 991--1004.

\bibitem{Area_et_al_2013} I. Area, D.K. Dimitrov, E. Godoy and V.G. Paschoa, Zeros of classical orthogonal polynomials of a discrete variable, {\it Math. Comp.} 82 (2013), 1069--1095.

\bibitem{Beardon}
A.F. Beardon, The theorems of Stieltjes and Favard, \emph{Comput. Methods Funct. Theory}  11(1) (2011), 247--262.

\bibitem{chiharabook} T.S. Chihara, \emph{An introduction to orthogonal polynomials}, Gordon and Breach, 1978.

\bibitem{Christoffel} E.B. Christoffel, \"Uber die Gaussische Quadratur und eine Verallgemeinerung derselben, \emph{J. Reine Angew. Math.} 55 (1858), 61--82.

\bibitem{DBS} C. de Boor and E.B. Saff, Finite sequences of orthogonal polynomials connected by a Jacobi matrix,
{\it Linear Algebra Appl.}  {75} (1986), 43--55.

\bibitem{Dimitrov}
D.K. Dimitrov, M.E.H. Ismail and F.R. Rafaeli,
Interlacing of zeros of orthogonal polynomials under modification of the measure, {\it J. Approx. Theory} 175 (2013), 64--76.

\bibitem{Dimitrov_Nikolov_2010} D.K. Dimitrov and  G.P. Nikolov, Sharp bounds for the extreme zeros of classical orthogonal polynomials, {\it J. Approx. Theory} 162 (2010), 1793--1804.

\bibitem{Dimitrov_Rafaeli_2009} D.K. Dimitrov and  F.R. Rafaeli, Monotonocity of zeros of Laguerre polynomials, {\it J. Comput. Appl. Math.} 233 (2009), 699--702.

\bibitem{DriverNumerMath} K. Driver, Interlacing of zeros of Gegenbauer polynomials of non-consecutive degree from different sequences. {\it Numer. Math.} 120 (2012), 35--44.

\bibitem{DJJ} K. Driver, A. Jooste and K. Jordaan. Stieltjes interlacing of zeros of Jacobi polynomials from different
sequences. {\it Electron. Trans. Numer. Anal.} 38 (2011), 317-326.

\bibitem{DJLag} K. Driver and K. Jordaan, Stieltjes interlacing of zeros of Laguerre polynomials from different sequences, {\it Indag. Math.} 21 (2011), 204--211.

\bibitem{DJ} K. Driver and K. Jordaan, Bounds for extreme zeros of some classical orthogonal polynomials, {\it J. Approx. Theory} 164 (2012), 1200--1204.

\bibitem{DM1} K. Driver and M.E. Muldoon, Zeros of pseudo-ultraspherical polynomials, {\it Anal. Appl.}, 12(5) (2012), 563--581.

\bibitem{Marcellan} J.C. Garcia-Ardila, F. Marcell\'an, P. H. Villamil-Hern\'andez, Associated orthogonal polynomials of the first kind and Darboux transformations. {\it J. Math. Anal. Appl.} 508(2) (2022):125883.

\bibitem{Gautschi} W. Gautschi,  An algorithmic implementation of the generalized Christoffel theorem, in
Numerical Integration (Hammerlin, G., ed.), ISNM 57, Birkhauser, Basel, 89--106.

\bibitem{Swami} P. Gochhayat, K. Jordaan, K Raghavendar and A. Swaminathan, Interlacing properties and bounds for zeros of $_2\Phi_1$ hypergeometric and little $q$-Jacobi polynomials, {\it Ramanujan J.} 40 (2016), 45--62.

\bibitem{Gupta_Muldoon_2007} D.E. Gupta and M.E. Muldoon, Inequalities for the smallest zeros of Laguerre polynomials and their $q$-analogues, {\it J. Inequal. Pure Appl. Math.} 8 (2007):24.

\bibitem{Ism} M.E.H. Ismail, {\it Classical and quantum orthogonal polynomials in one variable}, Encyclopedia of Mathematics and its Applications  98.  Cambridge: Cambridge University Press (2005).

\bibitem{Ismail_Li_1992} M.E.H. Ismail and X. Li,  Bound on the extreme zeros of orthogonal polynomials, {\it Proc. Amer. Math. Soc.} 115(1) (1992) 131--140.

\bibitem{ASJ} A.S. Jooste, A note on orthogonality and mixed recurrence equations, Constructive Theory of Functions, Sozopol 2019 (B. Draganov, K. Ivanov, G. Nikolov and R. Uluchev, Eds), 113--120. Prof. Marin Drinov Publishing House of BAS, Sofia, 2020.
\bibitem{Jooste_Jordaan} A. Jooste and K. Jordaan, Bounds for zeros of Meixner and Kravchuk polynomials, {\it LMS J. Comput. Math.} 17(1) (2014), 47--57.

\bibitem{JNK_2017} A. Jooste, P. Njionou Sadjang and W. Koepf, Inner bounds for the extreme zeros of $_3F_2$ hypergeometric polynomials, {\it Integral Transforms Spec. Funct.}  28(5) (2017), 361--373.

\bibitem{JT} K. Jordaan and F. To\'okos, Orthogonality and asymptotics of Pseudo-Jacobi polynomials for non-classical parameters, \emph{J.\ Approx.\ Theory} 178 (2014), 1--12.

\bibitem{Kar} P.P. Kar, K. Jordaan and P. Gochhayat, Steiltjes interlacing of zeros of little $q$-Jacobi and $q$-Laguerre polynomials from different sequences, {\it Numer. Algor.} 92 (2023), 723--746.

\bibitem{KLS} R. Koekoek, P.A. Lesky and R.F. Swarttouw, \emph{ Hypergeometric orthogonal polynomials and their $q$-analogues}. Springer Monogr. Math., Springer Verlag, Berlin (2010).

\bibitem{Krasikov_2002} I. Krasikov, Bounds for zeros of the Charlier polynomials, {\it Methods Appl. Anal.} 9(4) (2002), 599--610.
 
\bibitem{Krasikov_2006} I. Krasikov, On extreme zeros of classical orthogonal polynomials, {\it J. Comput. Appl. Math.} 193 (2006), 168--182.

\bibitem{Krasikov_Zarkh_2009} I. Krasikov and A. Zarkh, On zeros of discrete orthogonal polynomials, {\it J. Approx. Theory} 156 (2009), 121--141.

\bibitem{DLMF} F.W.J. Olver, A.B. Olde Daalhuis, D.W. Lozier, B.I. Schneider, R.F. Boisvert, C.W. Clark, B.R. Miller, B.V. Saunders, H.S. Cohl and M.A. McClain, eds. \emph{NIST Digital Library of Mathematical Functions.}, http://dlmf.nist.gov/, Release 1.1.7 of 2022-10-15.

\bibitem{rain} E.D. Rainville, \emph{Special Functions}, The Macmillan Company, New York, 1960.

\bibitem{szego} G. Szeg\H{o}, \emph{Orthogonal Polynomials}, American Mathematical Society, New York, 1959.

\bibitem{VZ} L. Vinet and A. Zhedanov, A characterization of classical and semiclassical orthogonal polynomials from their dual polynomials, {\it J. Comput. Appl. Math.} 172 (2004), 41-48.

\end{thebibliography}
\end{document}